\documentclass[11pt,twoside,a4paper]{article}

\usepackage[latin1]{inputenc}
\usepackage[T1]{fontenc}
\usepackage{a4wide}
\usepackage{amssymb,amsthm,amsmath,latexsym,stmaryrd}
\usepackage{graphicx}
\usepackage[right]{eurosym}
\usepackage[all]{xy}

\usepackage{pstricks}
\usepackage{pst-grad}
\usepackage{pst-text}
\usepackage{pst-tree}
\usepackage{pst-coil}
\usepackage{pst-node}

\title{On the Lie envelopping algebra of a pre-Lie algebra}

\author{J.-M. Oudom and D. Guin\\[2pt]
\small\hspace{7mm} oudom@math.univ-montp2.fr, dguin@math.univ-montp2.fr\\
\small I3M, Université Montpellier II, case 051\\
\small Place Eugène Bataillon\\
\small 34 095 Montpellier
}

\begin{document}

\theoremstyle{definition}
\newtheorem{definition}{Definition}[section]
\newtheorem{remarque}[definition]{Remark}
\newtheorem{notation}[definition]{Notation}
\newtheorem{exemple}[definition]{Example}
\newtheorem{df}{Definition}[subsection]
\newtheorem{ex}[df]{Example}
\newtheorem{rk}[df]{Remark}

\theoremstyle{plain}
\newtheorem{proposition}[definition]{Proposition}
\newtheorem{theorem}[definition]{Theorem}
\newtheorem{lemme}[definition]{Lemma}
\newtheorem{corollaire}[definition]{Corollary}
\newtheorem{prop}[df]{Proposition}
\newtheorem{thm}[df]{Theorem}

\def\dashededge{\ncline[linestyle=dashed,dash=2pt 2pt]}
\def\reddashededge{\ncline[linestyle=dashed,dash=2pt 2pt,linecolor=red]}

\maketitle

\quad For the last five years, many combinatorics Hopf algebras were introduced in different settings. One can quote the Hopf algebras of C. Brouder and A. Frabetti \cite{B-F}, A. Connes and D. Kreimer \cite{C-K}, L. Foissy \cite{F}, R. Grossman and R.-G. Larson \cite{G-L}, J.-L. Loday and M. Ronco  \cite{L-R}, I. Moerdijk and P. van der Laan\cite{Mo} \cite{VL}. 

\vspace{1mm}
\quad In this paper, our aim is to show that a lot of these Hopf algebras are related to  general algebraic constructions. In the commutative (or cocommutative) case, the key algebraic structure is the pre-Lie algebra one, as noticed by Chapoton and Livernet \cite{Ch-L}. In the non commutative and non cocommutative case, it seems to be played by brace algebras.

\vspace{1mm}
\quad In the first part, we recall the definition of a pre-Lie algebra and give some examples. The second part is devoted to the construction of an explicit $*$ product on the symmetric (co)algebra $S(L)$ of any pre-Lie algebra $L$. Then, we prove that $(S(L),*,\Delta)$ is isomorphic to the enveloping algebra of $L_{Lie}$. In the third part, we study the case where $L$ is the pre-Lie algebra of rooted trees, and we show that our construction corresponds to the dual of Connes and Kreimer's one. This gives another proof of the duality between the Connes-Kreimer and Grossman-Larson Hopf algebras \cite{C-K} \cite{P} \cite{Ho}. In the following part, we use the products $*$ and $\circ$ of the first part, in order to show that symmetric brace algebras introduced in \cite{L-M} are nothing else but pre-Lie algebras. Finally, in the fifth part, we recall the construction of a Hopf algebra structure on the tensor (co)algebra $T(V)$ of any brace algebra $V$ and we show that Foissy's  Hopf algebra of planar rooted trees belongs to this general setting. 

\vspace{1mm}
\quad From here we fix a commutative ring $\Bbbk$ over which modules, algebras, tensor products and linear maps are taken.

\section{Pre-Lie  algebras.}

\begin{definition}\quad

\vspace{1mm}
A pre-Lie  algebra is a module $L$ equipped with a bilinear product $\circ$ whose associator is symmetric in the two last variables :
\vspace{-2mm}
\begin{displaymath}
X\circ (Y\circ Z)-(X\circ Y)\circ Z\,=\,X\circ (Z\circ Y)-(X\circ Z)\circ Y.
\end{displaymath}
\end{definition}

These algebras were also called right symmetric algebras, Koszul-Vinberg or Vinberg algebras. Here, we use the terminology of F. Chapoton and M. Livernet in \cite{Ch-L}. There's no difficulty to give a graded version of pre-Lie algebras by replacing the above identity by :
 \vspace{-2mm}
\begin{displaymath}
X\circ (Y\circ Z)-(X\circ Y)\circ Z\,=\,(-1)^{|Y||Z|}(X\circ (Z\circ Y)-(X\circ Z)\circ Y).
\end{displaymath} 

\vspace{-2.5mm}
where the three variables are homogeneous and $|-|$ denotes the degree.

\vspace{3mm}
Let us now recall a well known fact, which is the root of the interest that people have shown in  pre-Lie algebras :

\begin{proposition}\quad

\vspace{1mm}
Let  $(L,\circ)$ be a pre-Lie algebra. The following bracket 
\vspace{-2mm}
\begin{displaymath}
[X,Y]=X\circ Y - Y\circ X
\end{displaymath}

\vspace{-2.5mm}
makes $L$ into a Lie algebra, which we will denote $L_{Lie}$.\qed 
\end{proposition}

In the following, we recall some classical and more recent examples.

\subsection{First examples.}

As its associator is symmetric in all three variables, an associative algebra is, of course,  a pre-Lie algebra. But the pre-Lie structure is in fact related to a weaker kind of associativity~: the associativity of the composition of multi-variables functions. M. Gerstenhaber already noted this fact in the Hochschild cohomology setting \cite{G} :

\begin{ex} Deformation complexes of algebras. 

\vspace{1mm}
Let $V$ be a module and let us denote $C^n(V,V)$ the space of all $n$-multilinear maps from $V$ to $V$. For $f\in C^p(V,V)$, $g\in C^q(V,V)$ and $i\in\llbracket 1,p\rrbracket$, one can define : 
\vspace{-2mm}
\begin{displaymath}
\begin{array}{l}f\circ_i g (x_1,\cdots\cdots, x_{p+q-1}) \,:=\, f(x_1,\cdots, x_{i-1},g(x_i,\cdots,x_{i+q-1}),x_{i+q},\cdots, x_{p+q-1})\\
f\circ g=\mathop{\sum}\limits_{i=1}^p f\circ_i g
\end{array}
\end{displaymath}

\vspace{-1mm}
Then the product  $\circ$ makes $C^\bullet(V,V)$ a pre-Lie algebra. In \cite{G}, M. Gerstenhaber gave a graded version of this product on the Hochschild cohomology complex  $C^\bullet(A,A)$ of an associative algebra $A$ :
\vspace{-2mm}
\begin{displaymath}
\begin{array}{l}f\circ_i g (x_1,\cdots, x_{p+q-1}) \,:=\, (-1)^{(p-1)(i-1)}f(x_1,\cdot\cdot, x_{i-1},g(x_i,\cdot\cdot,x_{i+q-1}),x_{i+q},\cdot\cdot, x_{p+q-1}).
\end{array}
\end{displaymath}

\vspace{-1mm}
Then, he showed that, in characteristic different from $2$, the graded Lie algebra induced by this graded pre-Lie structure controls the deformations of $A$.

More generally, D. Balavoine gave in  \cite{B} a similar graded Lie algebra construction  on the cohomology complex of an algebra over any quadratic operad, which controls  the given algebra's deformations.
\end{ex}

\begin{ex} Operads.

\vspace{1mm}
An operad  is a sequence of modules $(P(n)_{n\geq 2}$, where each $P(n)$ is a $k[\mathfrak{S_n}]$-module. For every $i$ between $1$ and $n$, we have a $\circ_i$ operation :
\vspace{-2mm}
\begin{displaymath}
\circ_i\,:\, P(n)\otimes P(m)\,\longrightarrow \,P(n+m-1)
\end{displaymath}

\vspace{-2mm}
satisfying the following conditions of associativity :
\vspace{-2mm}
\begin{displaymath}
\begin{array}{c}
f\circ_i (g\circ_j h)\, = \, (f\circ_i g)\circ_{i+j-1}h\\
(f\circ_i g) \circ_{|g|+i+j-1} h\,=\, (f\circ_j h) \circ_{|h|+i+j-1} g
\end{array}
\end{displaymath}

\vspace{-2mm}
This $\circ_i$-operations should moreover be some way compatible with the $k[\mathfrak{S_n}]$-actions. There is no difficulty to check that the product :
\vspace{-2mm}
\begin{displaymath}
f\circ g\,=\, \sum_{i=1}^{|f|}f\circ_i g
\end{displaymath}

\vspace{-2mm}
defines a pre-Lie product on $\mathop{\bigoplus}\limits_{n\geq 2} P(n)$.
\vspace{1mm}
\end{ex}

\subsection{A geometric example : affine manifolds.}

An affine manifold is a manifold with a torsion free and flat connection $\nabla$. Equivalently, it is a manifold equipped with an atlas with affine transition functions. The vector fields of such a manifold is a pre-Lie algebra for the following circle product :

\vspace{-3mm}
$$
X\circ Y = \nabla_X(Y),
$$

\vspace{-2mm}
whose associated Lie bracket is the usual one.

\subsection{The inspiring example : rooted trees.}

\begin{df}
 
\vspace{1mm}
 A rooted tree is a tree with a distinguished vertex : its root.  A combinatorial definition could be the following : it is a finite poset with a minimum (the root) and no critical pair :
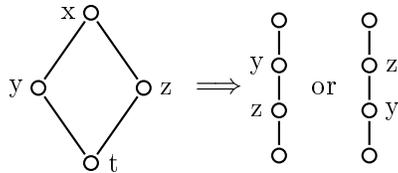
\vspace{-2mm}

\rput(4.7,0.1){x}
\rput(4,-0.9){y}
\rput(6,-0.9){z}
\rput(5.3,-1.9){t}
\rput(7.2,-0.6){y}
\rput(7.2,-1.2){z}
\rput(9,-0.6){z}
\rput(9,-1.2){y}
\cnode(5,0.1){3pt}{A}
\cnode(4.3,-0.9){3pt}{B}
\cnode(5.7,-0.9){3pt}{C}
\cnode(5,-1.9){3pt}{D}
\ncline[nodesep=1pt]{A}{B}
\ncline[nodesep=1pt]{A}{C}
\ncline[nodesep=1pt]{D}{B}
\ncline[nodesep=1pt]{D}{C}
\rput(6.7,-0.9){$\Longrightarrow$}
\cnode(7.5,-0.6){3pt}{E}
\cnode(7.5,-1.2){3pt}{F}
\cnode(7.5,0){3pt}{E'}
\cnode(7.5,-1.8){3pt}{F'}
\ncline[nodesep=1pt]{E}{F}
\ncline[nodesep=1pt]{E}{E'}
\ncline[nodesep=1pt]{F}{F'}
\rput(8.1,-0.9){or}
\cnode(8.7,-0.6){3pt}{G}
\cnode(8.7,-1.2){3pt}{H}
\cnode(8.7,0){3pt}{G'}
\cnode(8.7,-1.8){3pt}{H'}
\ncline[nodesep=1pt]{G}{H}
\ncline[nodesep=1pt]{G}{G'}
\ncline[nodesep=1pt]{H}{H'}
\vspace{2cm}

For any rooted tree $T$, we denote $|T|$ the underlying set and we call its vertices of $T$ its elements.  Let $X$ be a set. A $X$-colored rooted tree is a rooted tree $T$ equipped with a \emph{color} map $|T|\rightarrow X$, which associates its color in $X$ to every vertex of $T$. 
\end{df}

\vspace{1mm}
We will represent trees by using planar graphs and the order induced by the gravity : 

\vspace{1mm}
\rput(5.5,1){\Large $=$}
\rput(8.9,1){\Large $=$}
\hspace{2.5cm}\pstree[treemode=U,treesep=5mm,levelsep=1cm]{\Tc{3pt}}{
 \pstree{\Tc{3pt}}{\pstree{\Tc{3pt}}{\Tc{3pt} \Tn} \Tc{3pt}}
\pstree{\Tc{3pt}}{\Tn\Tc{3pt}}}
\,\qquad \,
\pstree[treemode=U,treesep=5mm,levelsep=1cm]{\Tc{3pt}}{
\pstree{\Tc{3pt}}{\Tc{3pt} \Tn} \pstree{\Tc{3pt}}{\Tc{3pt} \pstree{\Tc{3pt}}{\Tn \Tc{3pt}}}}
\,\qquad \,
\pstree[treemode=U,treesep=5mm,levelsep=1cm]{\Tc{3pt}}{
\pstree{\Tc{3pt}}{\Tc{3pt} \Tn} \pstree{\Tc{3pt}}{\pstree{\Tc{3pt}}{\Tc{3pt} \Tn}\Tc{3pt}}}

\vspace{1mm}
Here are two different colored trees with the same underlying tree :

\vspace{1mm}
\rput(7.2,1){\Large $\neq$}
\hspace{3.5cm}
\pstree[treemode=U,treesep=5mm,levelsep=1cm]{\Tc*[linecolor=gray]{3pt}}{\pstree{\Tc*[linecolor=gray]{3pt}}{\Tc{3pt} \Tc*[linecolor=black]{3pt}}\pstree{\Tc*[linecolor=black]{3pt}}{\Tc*[linecolor=gray]{3pt}} \pstree{\Tc{3pt}}{\Tc{3pt} \Tc*[linecolor=black]{3pt}}}
\,\qquad \,
\pstree[treemode=U,treesep=5mm,levelsep=1cm]{\Tc*[linecolor=black]{3pt}}{\pstree{\Tc*[linecolor=black]{3pt}}{\Tc{3pt} \Tc*[linecolor=gray]{3pt}}\pstree{\Tc*[linecolor=black]{3pt}}{\Tc{3pt}} \pstree{\Tc*[linecolor=gray]{3pt}}{\Tc*[linecolor=black]{3pt} \Tc{3pt}}}

\vspace{2mm}
For two given  rooted trees $T_1$ and $T_2$, and a chosen vertex $v$ of $T_1$, we can glue $T_1$ on $v$. The rooted tree $T_1\circ_v T_2$ obtained in such a way is the poset whose underlying set is the disjoint union  $|T_1|\amalg |T_2|$. The order is induced by the orders of $T_1$ and $T_2$ and $w>v$ for all $w$ in $T_2$ :

\vspace{0.9cm}
\rput(3.1,1){$\scriptstyle v$}
\rput(9.2,0){$\scriptstyle v$}
\rput(8.2,1){$:=$}
\rput(5.5,0.5){$\circ_v$}
\hspace{3cm}\pstree[treemode=U,treesep=5mm,levelsep=1cm]{\Tc{3pt}}{
 \pstree{\Tc{3pt}}{\Tc{3pt} \Tn } \Tc{3pt} \Tc{3pt}}
\,\qquad \,
\pstree[treemode=U,treesep=5mm,levelsep=1cm]{\Tc{3pt}}{
\Tc{3pt} \Tn \pstree{\Tc{3pt}}{\Tc{3pt} \Tc{3pt}}}
\,\qquad \,
\rput(1,1){\pstree[treemode=U,treesep=5mm,levelsep=1cm]{\Tc{3pt}}{
 \pstree{\Tc{3pt}}{\Tc{3pt} \Tn \pstree{\Tc[edge=\dashededge]{3pt}}{\Tc{3pt} \Tn \pstree{\Tc{3pt}}{\Tc{3pt} \Tc{3pt}}}} \Tc{3pt} \Tc{3pt}}}

\vspace{1.2cm}
Notice that, when $T_1$ and $T_2$ are colored by a set $X$, then $T_1\circ_v T_2$ is a $X$-colored rooted tree, whose \emph{color} map is the disjoint union of the two \emph{color} maps of $T_1$ and $T_2$.

\begin{prop}\quad

\vspace{1mm}
Let us denote $\mathcal P\mathcal L (X)$ the free module whose base is the set of $X$-colored rooted trees. We define on $\mathcal P\mathcal L (X)$ the following $\circ$ product :
\vspace{-2mm}
\begin{displaymath}
T_1\circ T_2\,=\, \sum_{v\in|T_1|} T_1\circ_v T_2
\end{displaymath}

\vspace{-6mm}
Then $(\mathcal P\mathcal L (X),\circ)$ is a pre-Lie algebra.
\end{prop}

\begin{proof} It follows from the following identity : 
\vspace{-2mm}
\begin{displaymath}
(T_1\circ T_2)\circ T_3 - T_1\circ (T_2\circ T_3) \,=\, \sum_{v,w\in |T_1|} (T_1\circ_v T_2)\circ_w T_3 \, = \, \sum_{v,w\in |T_1|} (T_1\circ_w T_3)\circ_v T_2
\end{displaymath}
\vspace{-4mm}

which can be illustrated by the the following picture :

\pnode(2.2,1){A}
\pnode(2.2,0){B}
\ncarc[arcangleA=45,arcangleB=45]{-}{B}{A}
\pnode(4.9,1){C}
\pnode(4.9,0){D}
\ncarc[arcangleA=45,arcangleB=45]{-}{C}{D}
\pnode(8.45,1){E}
\pnode(8.45,0){F}
\ncarc[arcangleA=45,arcangleB=45]{-}{F}{E}
\pnode(11.6,1){G}
\pnode(11.6,0){H}
\ncarc[arcangleA=45,arcangleB=45]{-}{G}{H}
\rput(3.5,0.5){\Large $\circ$}
\rput(5.3,0.5){\Large $\circ$}
\rput(6.6,0.5){\Large $-$}
\rput(8,0.5){\Large $\circ$}
\rput(10,0.5){\Large $\circ$}
\hspace{2.2cm}\;\pstree[treemode=U,treesep=5mm,levelsep=1cm]{\Tc{3pt}}{\Tc{3pt} \Tc{3pt}}  \; \pstree[treemode=U,treesep=5mm,levelsep=1cm]{\Tc{3pt}}{\Tdia{$T_1$}} \;\;\;  \pstree[treemode=U,treesep=5mm,levelsep=1cm]{\Tc{3pt}}{\Tdia{$T_2$}} \;\; \pstree[treemode=U,treesep=5mm,levelsep=1cm]{\Tc{3pt}}{\Tc{3pt} \Tc{3pt}} \;\;  \;\;\,\pstree[treemode=U,treesep=5mm,levelsep=1cm]{\Tc{3pt}}{\Tdia{$T_1$}} \;\;  \pstree[treemode=U,treesep=5mm,levelsep=1cm]{\Tc{3pt}}{\Tdia{$T_2$}}

\vspace{2mm}
\rput(0.2,1){\Large $=$}
\rput(3.9,1.1){ $+\;\;2\times $}
\rput(6.6,1.1){ $+\;\;2\times $}
\rput(9.6,1.1){ $+\;\;2\times $}
\rput(12.8,1.1){ $+\;\;2\times $}
\hspace{2mm} \pstree[treemode=U,treesep=5mm,levelsep=1cm]{\Tc{3pt}}{\pstree{\Tc[edge=\dashededge]{3pt}}{\Tdia{$T_1$}}\Tc{3pt}\Tc{3pt}\pstree{\Tc[edge=\dashededge]{3pt}}{\Tdia{$T_2$}}}
\hspace{-1.5cm}\pstree[treemode=U,treesep=5mm,levelsep=1cm]{\Tc{3pt}}{\pstree{\Tc{3pt}}{\Tn\Tn\Tn\Tn\pstree{\Tc[edge=\dashededge]{3pt}}{\Tdia{$T_1$}}\pstree{\Tc[edge=\dashededge]{3pt}}{\Tdia{$T_2$}}\Tn\Tn\Tn\Tn} \Tc{3pt}}
\hspace{-1.5cm}\pstree[treemode=U,treesep=5mm,levelsep=1cm]{\Tc{3pt}}{\Tn\Tn\pstree{\Tc{3pt}}{ \pstree{\Tc[edge=\dashededge]{3pt}}{\Tdia{$T_1$}}\Tn} \pstree{\Tc{3pt}}{\Tn\pstree{\Tc[edge=\dashededge]{3pt}}{\Tdia{$T_2$}}}\Tn\Tn}
\pstree[treemode=U,treesep=5mm,levelsep=1cm]{\Tc{3pt}}{\pstree{\Tc{3pt}}{\pstree{\Tc[edge=\dashededge]{3pt}}{\Tdia{$T_1$}}}\Tc{3pt}\pstree{\Tc[edge=\dashededge]{3pt}}{\Tdia{$T_2$}}}
\hspace{5mm}\pstree[treemode=U,treesep=5mm,levelsep=1cm]{\Tc{3pt}}{\pstree{\Tc[edge=\dashededge]{3pt}}{\Tdia{$T_1$}}\Tc{3pt}\pstree{\Tc{3pt}}{\pstree{\Tc[edge=\dashededge]{3pt}}{\Tdia{$T_2$}}}}
\end{proof}

\begin{rk}  F. Chapoton and M. Livernet have shown in \cite{Ch-L} that $\mathcal{PL}(X)$ is the free pre-Lie algebra generated by $X$. We will give a (not so) different proof of this fact in section 3.
\end{rk}
\subsection{Other graphical examples.}

As pointed out in the first examples, the pre-Lie algebra structure is related to a kind of \emph{weak} associativity, a symmetric substitution associativity. This \emph{weak} associativity often appears in the world of graphs. The previous example is based on graftings, and we can extend this grafting product to graphs with a marked vertex. For two pointed graphs $(G_1,\bullet_1)$ and $(G_2,\bullet_2)$, one can graft $\bullet_2$ on each vertex of $G_1$. If we take the sum of all this grapftings, we get a pre-Lie product on the free modules generated by the finite pointed graphs :

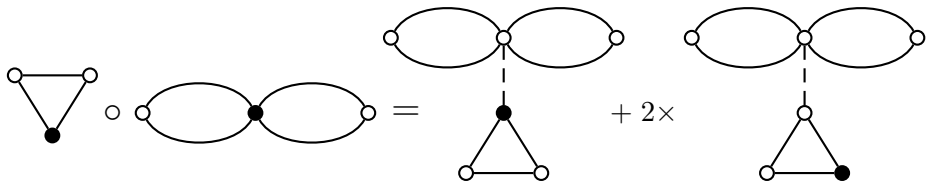
\cnode(3.2,-1.3){3pt}{A}
\cnode*(4.7,-1.3){3pt}{B}
\cnode(6.2,-1.3){3pt}{C}
\cnode(1.5,-0.8){3pt}{E}
\cnode(2.5,-0.8){3pt}{F}
\cnode*(2,-1.6){3pt}{G}
\cnode(6.5,-0.3){3pt}{a}
\cnode(8,-0.3){3pt}{b}
\cnode(9.5,-0.3){3pt}{c}
\cnode(7.5,-2.1){3pt}{e}
\cnode(8.5,-2.1){3pt}{f}
\cnode*(8,-1.3){3pt}{g}
\cnode(10.5,-0.3){3pt}{aa}
\cnode(12,-0.3){3pt}{bb}
\cnode(13.5,-0.3){3pt}{cc}
\cnode(11.5,-2.1){3pt}{ee}
\cnode*(12.5,-2.1){3pt}{ff}
\cnode(12,-1.3){3pt}{gg}
\rput(2.8,-1.3){\Large $\circ$}
\rput(6.7,-1.3){\Large $=$}
\rput(9.8,-1.3){ $+\; 2\times$}
\ncarc[arcangleA=60,arcangleB=60]{-}{A}{B}
\ncarc[arcangleA=60,arcangleB=60]{-}{B}{A}
\ncarc[arcangleA=60,arcangleB=60]{-}{C}{B}
\ncarc[arcangleA=60,arcangleB=60]{-}{B}{C}
\ncarc[arcangleA=60,arcangleB=60]{-}{a}{b}
\ncarc[arcangleA=60,arcangleB=60]{-}{b}{a}
\ncarc[arcangleA=60,arcangleB=60]{-}{c}{b}
\ncarc[arcangleA=60,arcangleB=60]{-}{b}{c}
\ncarc[arcangleA=60,arcangleB=60]{-}{aa}{bb}
\ncarc[arcangleA=60,arcangleB=60]{-}{bb}{aa}
\ncarc[arcangleA=60,arcangleB=60]{-}{cc}{bb}
\ncarc[arcangleA=60,arcangleB=60]{-}{bb}{cc}
\ncline{E}{F}
\ncline{E}{G}
\ncline{G}{F}
\ncline{e}{f}
\ncline{e}{g}
\ncline{g}{f}
\ncline{ee}{ff}
\ncline{ee}{gg}
\ncline{gg}{ff}
\ncline[linestyle=dashed]{b}{g}
\ncline[linestyle=dashed]{bb}{gg}

\vspace{2.5cm}
Moreover, one can form pre-Lie algebra structures on graphs by using an insertion process.  For instance, one can define a pre-Lie algebra product on the free module generated by the finite graphs, where the product of two graphs $G_1$ and $G_2$ is the sum of all possible insertions of $G_2$ in $G_1$. Here, an insertion of $G_2$ in $G_1$ consists in choosing a non univalent vertex $v$ of $G_1$ and a bijection $\varphi$ from the set of univalent vertices of $G_2$ to the set of edges of $v$. To obtain the resulting graph, we suppress the vertex $v$ and glue the univalent vertices to the edges of $v$ according to $\varphi$ :

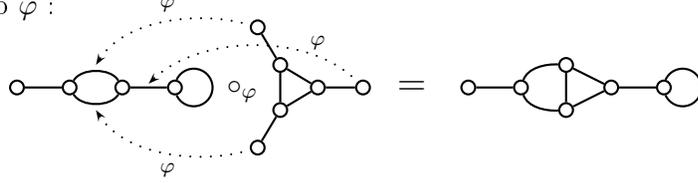
\rput(5.5,-0.55){$\circ_\varphi$}
\rput(6.5,0.1){$\scriptstyle \varphi$}
\rput(4.5,0.6){$\scriptstyle \varphi$}
\rput(4.5,-1.6){$\scriptstyle \varphi$}
\rput(7.75,-0.5){\Large $=$}
\cnode(2.5,-0.5){3pt}{A}
\cnode(3.2,-0.5){3pt}{B}
\cnode(3.9,-0.5){3pt}{C}
\cnode(4.6,-0.5){3pt}{D}
\pnode(3.55,-0.2){E}
\pnode(3.55,-0.8){F}
\pnode(4.25,-0.45){G}
\cnode(6,-0.2){3pt}{a}
\cnode(6,-0.8){3pt}{b}
\cnode(6.5,-0.5){3pt}{c}
\cnode(5.7,0.3){3pt}{d}
\cnode(5.7,-1.3){3pt}{e}
\cnode(7.1,-0.5){3pt}{f}
\cnode(8.5,-0.5){3pt}{aa}
\cnode(9.2,-0.5){3pt}{bb}
\cnode(9.8,-0.2){3pt}{cc}
\cnode(9.8,-0.8){3pt}{dd}
\cnode(10.4,-0.5){3pt}{ee}
\cnode(11.1,-0.5){3pt}{ff}
\ncline{A}{B}
\ncarc[arcangleA=60,arcangleB=60]{-}{C}{B}
\ncarc[arcangleA=60,arcangleB=60]{-}{B}{C}
\ncline{C}{D}
\nccircle[angleA=-90]{D}{0.25}
\ncline{a}{b}
\ncline{a}{c}
\ncline{c}{b}
\ncline{a}{d}
\ncline{e}{b}
\ncline{c}{f}
\ncarc[arcangleA=45,arcangleB=30,linestyle=dotted]{<-}{E}{d}
\ncarc[arcangleA=30,arcangleB=45,linestyle=dotted]{->}{e}{F}
\ncarc[arcangleA=45,arcangleB=45,linestyle=dotted]{<-}{G}{f}
\nccircle[angleA=-90]{ff}{0.25}
\ncline{aa}{bb}
\ncline{ee}{ff}
\ncline{cc}{dd}
\ncarc[arcangleA=30,arcangleB=45]{-}{dd}{bb}
\ncarc[arcangleA=45,arcangleB=30]{-}{bb}{cc}
\ncline{dd}{ee}
\ncline{cc}{ee}

\vspace{1.7cm}
F. Chapoton has proposed a similar insertion example in \cite{Ch1}.

\section{The $\boldsymbol{\circ}$ and  $\boldsymbol{*}$ products on  $\boldsymbol{S(L)}$.}

\begin{notation}\quad

\vspace{1mm}
Let $L$ be a pre-Lie algebra. We denote $S(L)$ its symmetric algebra with the usual shuffle coproduct denoted $\Delta$. We will use without restraint the classical  Sweedler's notations : $\Delta (A)=A_{(1)}\otimes A_{(2)}$.
\end{notation}

Our aim in the following is to extend the $\circ$ product of $L$ to the whole $S(L)$. Then,  we will be able to define an associative product $*$ on $S(L)$. First, let us extend $\circ$ as a right $L$-action on $S(L)$ :  

\begin{notation}\quad

\vspace{1mm}
Let  $T$, $X_1$, \dots, $X_n$ be in $L$, and set :
\vspace{-2mm}
\begin{displaymath}
1\circ T := 0 \quad \hbox{\rm and }\quad  (X_1\cdots X_n)\circ T\,:=\,\sum_{1\leq i\leq n} \cdots ( X_i\circ T)\cdots
\end{displaymath}

\vspace{-3mm}
\end{notation}

\begin{remarque} It's not difficult to see that this action is right symmetric : 
\vspace{-2mm}
\begin{displaymath}
A\circ (X\circ Y) - (A\circ X)\circ Y \,=\, A\circ (Y\circ X) - (A\circ Y)\circ X 
\end{displaymath}

\vspace{-3mm}
where $A\in S(L)$ and $X$, $Y$ are in $L$.

\end{remarque}
\vspace{1mm}

Now, we define $T\circ A$ where $A$ is a monomial of $S(L)$ and $T$ belongs to $L$ :

\begin{definition}\quad

\vspace{1mm}
We define inductively multilinear maps :
\vspace{-2mm}
\begin{displaymath}
\begin{array}{ccccc}
\circ_{\scriptscriptstyle (n)}&:& L\otimes L^{\otimes n} & \longrightarrow & L\\
&& T\otimes A & \longmapsto & T\circ A
\end{array}
\end{displaymath}

by the following way : $T\circ 1 \,:=\, T$ and for $X\in L$, $B$ in $S(L)$ and $A=BX$,
\vspace{-2mm}
\begin{displaymath}
T\circ_{\scriptscriptstyle (n)} A \,:=(T\circ_{\scriptscriptstyle (n-1)}B)\circ X -T\circ_{\scriptscriptstyle (n-1)} (B\circ X).
\end{displaymath}

\vspace{-3mm}
\end{definition}

\begin{lemme}\quad

\vspace{1mm}
For every integer $n$, the map $\circ_{\scriptscriptstyle (n)}$ is symmetric in the $n$ last variables. Therefore, it defines a map : 
\vspace{-2mm}
\begin{displaymath}
\begin{array}{ccccc}
\circ_{\scriptscriptstyle (n)}&:& L\otimes S^nL & \longrightarrow & L\\
&& T\otimes A & \longmapsto & T\circ A
\end{array}
\end{displaymath}
\end{lemme}

\begin{proof}
The invariance by permutation of the last variables is obtained by induction on $n$ : 
\vspace{-2mm}
\begin{displaymath}
\begin{array}{ccl}
T\circ_{(n+2)} AXY &=& (T\circ_{(n+1)} AX)\circ Y \,-\, T\circ_{(n+1)} (AX\circ Y)\\
&=& ((T\circ_{(n)} A)\circ X)\circ Y - (T\circ_{(n)} (A\circ X))\circ Y \\
& &- T\circ_{(n+1)} (A\circ Y)X - T\circ_{(n+1)} A(X\circ Y)\\
&=&  ((T\circ_{(n)} A)\circ X)\circ Y - (T\circ_{(n)} (A\circ X))\circ Y\\
&& - (T\circ_{(n)} (A\circ Y))\circ X + T\circ_{(n)} ((A\circ Y)\circ X)\\
&& -(T\circ_{(n)} A)\circ (X\circ Y) + T\circ_{(n)} (A\circ (X\circ Y))
\end{array}
\end{displaymath}

In the above equality, it follows from the pre-Lie identity that the first terms of the first and third lines give rise to an expression symmetric in $X$ and $Y$. By (2.3), the sum of the last terms of the two last lines is also symmetric in $X$ and $Y$. Finally, the sum of the two remaining terms is obviously symmetric in $X$ and $Y$.

By this way, we get the invariance by the transposition $(n-1,n)$. If we assume the symmetry of $\circ_{(n-1)
}$, we have the symmetry in the $n-1$ variables of $AX$. Thus, as  $\mathfrak{S}_n$ is generated by $(n-1,n)$ and $\mathfrak{S}_{n-1}\times \mathfrak{S}_{1}$, we get the announced invariance. 
\end{proof}

\begin{remarque} In terms of operads, these $\circ_{(n)}$ operations were introduced by F. Chapoton and M. Livernet in \cite{Ch-L}.
\end{remarque}

\begin{proposition}\quad

\vspace{1mm}
There's a unique extension of the  product $\circ$ to $S(L)$  such that :
\vspace{-2mm}
\begin{displaymath}\left\{\begin{array}{cl}
(i)   & A\circ 1\, = \, A\\
(ii)  & T\circ BX\, = \, (T\circ B)\circ X\, -\,T\circ (B\circ X) \\
(iii) & AB\circ C \, = \, (A\circ C_{(1)})(B\circ C_{(2)})
\end{array}\right.
\end{displaymath}

\vspace{-2mm}
where  $A$, $B$, $C$ belong to  $S(L)$ and $X$ is in $L$.
\end{proposition}

\begin{proof} First, notice that, according to $(i$) and $(iii)$, we necessarily have   $1\circ T = 0$ for every $T$ in $L$. Thus, we have, by induction and $(iii)$,  $1\circ A = \varepsilon (A)$ for all $A$ in $S(L)$, where $\varepsilon$ denotes the counit of $(S(L),\Delta)$. 

We also get by induction,  $(i)$ and  $(iii$)  : 
\vspace{-2mm}
\begin{displaymath}
(X_1\cdots X_n)\circ T\,=\,\sum_{1\leq i\leq n} \cdots ( X_i\circ T)\cdots
\end{displaymath}

\vspace{-3mm}
It follows from $(ii)$ that, for all monomial $A$ of length $n$ : $T\circ A \,=\,T\circ_{\scriptscriptstyle (n)} A$.

Then, $(iii)$ gives the definition of $A\circ B$ for any monomials $A$ and $B$. Thanks to the coassociativity and cocommutativity of $\Delta$, this definition is not ambiguous : 
\vspace{-2mm}
\begin{displaymath}
\begin{array}{ccl}
(A\circ C_{(1)})(B\circ C_{(2)}) &=& (B\circ C_{(2)})(A\circ C_{(1)})\\
&=& (B\circ C_{(1)})(A\circ C_{(2)})\\

(AB\circ D_{(1)})(C\circ D_{(2)}) &=& (A\circ D_{(1)})(B\circ D_{(2)})(C\circ D_{(3)})\\
&=& (A\circ D_{(1)})(BC\circ D_{(2)})
\end{array}
\end{displaymath}
\end{proof}

\begin{proposition}\quad

\vspace{1mm}
Let $A$, $B$, $C$  be in $S(L)$ and let  $X$ be in  $L$. We have :
\vspace{-2mm}
\begin{displaymath}
\begin{array}{cl} 
(i) & 1\circ A \, = \, \varepsilon (A)\\
(ii) & \varepsilon (A\circ B)\,=\, \varepsilon (A)\varepsilon (B)\\
(iii) & \Delta(A\circ B)\, = \, (A_{(1)}\circ B_{(1)}) \otimes (A_{(2)}\circ B_{(2)})\\
(iv)   & A\circ BX \, = \, (A\circ B)\circ X - A\circ (B\circ X),\\
(v)  & (A\circ B)\circ C\, =\, A\circ ((B\circ C_{(1)})C_{(2)}). 
\end{array}
\end{displaymath}
\end{proposition}

\begin{proof} We have already seen $(i)$ in the proof of (2.6.). The $(ii)$ part follows from the fact that $1\circ 1 = 1$. Moreover,  when $A$ is a monomial, $A\circ B$ has the same length as $A$.

\vspace{1mm}
We  get  $(iii)$ by induction on the  length of $A$.  When $A$ is of length $0$, we deduce from (i) that~: $\Delta (A\circ B)\,=\, \varepsilon (A)\varepsilon (B) 1 \otimes 1  \, = \Delta (A)\circ \Delta(B)$. Then, we use (2.6.$iii$) to decrease the length of $A$ : 
\vspace{-2mm}
\begin{displaymath}
\begin{array}{ccll} 
\Delta(AB\circ C) & = & \Delta((A\circ C_{(1)})(B\circ C_{(2)}))\\
&=&\Delta(A\circ C_{(1)})\Delta(B\circ C_{(2)})&\\
&=& (A_{(1)}\circ C_{(1)}\otimes A_{(2)}\circ C_{(2)}) (B_{(1)}\circ C_{(3)}\otimes B_{(2)}\circ C_{(4)})\\
&=& (A_{(1)}\circ C_{(1)})(B_{(1)}\circ C_{(2)})\otimes  (A_{(2)}\circ C_{(3)})(B_{(2)}\circ C_{(4)})\\
&=& (A_{(1)} B_{(1)}\circ C_{(1)})\otimes (A_{(2)}B_{(2)}\circ C_{(2)})\qquad\hbox{\rm from (2.6.}iii)\\
&=& \Delta(AB)\circ \Delta(C).&
\end{array}
\end{displaymath}

We get $(iv)$ from (2.6.$ii$),  (2.6.$iii$) and the above $(iii)$ by induction on the length of $A$. When $A$ is of length $0$, everything is null. 

\vspace{-2mm}
\begin{displaymath}
\begin{array}{ccl} 
AB\circ CX & = & (A\circ C_{(1)}X) (B\circ C_{(2)}) + (A\circ C_{(1)})(B\circ C_{(2)}X) \\
& = &  ((A\circ C_{(1)})\circ X)(B\circ C_{(2)})   - (A\circ (C_{(1)}\circ X))(B\circ C_{(2)}) \\
&& + (A\circ C_{(1)})  ((B\circ C_{(2)})\circ X)  - (A\circ C_{(1)})  (B\circ (C_{(2)}\circ X))\\
&=& ((A\circ C_{(1)})(B\circ C_{(2)}))\circ X - (A\circ (C_{(1)}\circ X))(B\circ C_{(2)})\\
&& - (A\circ C_{(1)})  (B\circ (C_{(2)}\circ X))\\
&=& (AB\circ C)\circ X - AB\circ (C\circ X).\\
\end{array}
\end{displaymath}

\vspace{-2mm}
Finally, $(v)$ also follows from an induction on the length of $A$. The $0$ length case follows from $(i)$ and $(ii)$. Then,  following (2.6.$iii$), 
\vspace{-2mm}
\begin{displaymath}
\begin{array}{ccl} 
(AB\circ C)\circ D & = & ((A\circ C_{(1)})\circ D_{(1)}) ((B \circ C_{(2)})\circ D_{(2)}) \\
& = & (A\circ ((C_{(1)}\circ D_{(1)})D_{(2)}))(B \circ ((C_{(2)}\circ D_{(3)})D_{(4)}))\\
& = & AB\circ ((C \circ D_{(1)})D_{(2)}),
\end{array}
\end{displaymath}

where the last equality follows from (2.6.$iii$) and :
\vspace{-2mm}
\begin{displaymath}
\begin{array}{ccl} 
\Delta ((C \circ D_{(1)})D_{(2)}) &= & (\Delta (C)\circ \Delta (D_{(1)})) \Delta (D_{(2)})\\
& = & (C_{(1)}\circ D_{(1)})D_{(3)} \otimes (C_{(2)}\circ D_{(2)}) D_{(4)}\\
& = & (C_{(1)}\circ D_{(1)})D_{(2)} \otimes (C_{(2)}\circ D_{(3)}) D_{(4)}
\end{array}
\end{displaymath}
\end{proof}

\begin{definition}\quad

\vspace{1mm}
We define the following $*$ product  on $S(L)$ :
\vspace{-2mm}
\begin{displaymath}
A*B\,:=\, (A\circ B_{(1)})B_{(2)}.
\end{displaymath}
\end{definition}

\begin{lemme}\quad

\vspace{1mm}
The product $*$ is associative and makes $S(L)$ into a Hopf algebra with coproduct  $\Delta$.
\end{lemme}

\begin{proof}  The compatibility between $*$ and $\Delta$ follows from (2.7.$iii$). Let us prove the associativity of $*$ :

\vspace{-2mm}
\begin{displaymath}
\begin{array}{ccll} 
(A*B)*C & = & (((A\circ B_{(1)})B_{(2)})\circ C_{(1)})C_{(2)} &\hspace{-5mm}\hbox{\rm by definition of } *, \\\
& = & ((A\circ B_{(1)})\circ C_{(1)})(B_{(2)}\circ C_{(2)})C_{(3)} &\hspace{-3mm}\hbox{\rm from (2.6.}iii),\\
&=& (A\circ ((B_{(1)})\circ C_{(1)})C_{(2)}))(B_{(2)}\circ C_{(3)})C_{(4)}&\hbox{\rm from (2.7.}v),\\
&=&(A\circ ((B_{(1)})\circ C_{(1)})C_{(3)}))(B_{(2)}\circ C_{(2)})C_{(4)}&\hbox{\rm by cocommutivity of } \Delta,\\
&=& A*((B\circ C_{(1)})C_{(2)})&\hspace{-15mm}\hbox{\rm by definition of }*~\hbox{\rm and (2.7.}iii),\\
&=& A*(B*C)&
\end{array}
\end{displaymath}
\end{proof}

\begin{lemme}\quad

\vspace{1mm}
The product  $*$ agrees with the usual increasing filtration  of $S(L)$ by the length of the words, and the associated graded product is the usual product of the symmetric algebra. 
\end{lemme}

\begin{proof} The first part follows from the fact that the length of  $A\circ B$ is the same as the length of $A$.  For the second part, we check that, for any monomials $A$ and $B$, the highest length term of $A*B$ is $AB$.
\end{proof}

\begin{theorem}\quad

\vspace{1mm}
The Hopf algebra $(S(L), *,\Delta)$ is isomorphic to the enveloping algebra  $\mathcal U(L_{Lie})$.
\end{theorem}

\begin{proof} For any $X$ and $Y$ in $L$, we have :
\vspace{-3mm}
\begin{displaymath}
\begin{array}{ccl} 
X*Y-Y*X&=& XY+X\circ Y - YX -Y\circ X\\
&=& [X,Y]
\end{array}
\end{displaymath}

\vspace{-3mm}
Therefore, the morphism from $T(L)$ to $S(L)$ induced by the inclusion of $L$ factors out through the ideal $<X\otimes Y - Y\otimes X -[X,Y]>$ :
\vspace{-3mm}
$$\xymatrix{
(T(L),\otimes,\Delta)\ar[rr]\ar[dr]_\pi &  &(S(L),*,\Delta)\\
&(\mathcal U(L_{Lie}),\cdot,\Delta)\ar[ur]_\varphi&
}
$$

\vspace{-3mm}
By (2.10.) all the morphisms involved are filtered by the length of the words. Moreover, as $T(L)$ and $\mathcal U(L_{Lie})$ are generated by $L$, they are Hopf algebra's morphisms. By taking the corresponding graded maps, we obtain the following diagramm of Hopf algebras :
\vspace{-3mm}
$$\xymatrix{
T(L)\ar[rr]\ar[dr]_{\overline{\pi}} &  &S(L)\\
&gr U(L_{Lie})\ar[ur]_{\overline{\varphi}}&
}
$$

\vspace{-3mm}
The horizontal map is nothing but the usual projection of $T(L)$ on $S(L)$. Moreover, $\overline{\pi}$ factors out through $S(L)$, which gives :
\vspace{-2mm}
$$\xymatrix{
S(L)\ar[rr]^{id}\ar[dr]_{\overline{\pi}} &  &S(L)\\
&gr U(L_{Lie})\ar[ur]_{\overline{\varphi}}&
}
$$
\vspace{-3mm}

As all the maps involved are onto, $\overline{\varphi}$ is an isomorphism. Therefore, $\varphi$ is an isomorphism.
\end{proof}

\begin{remarque} Everything we have done has naturally a graded version. In fact, all can be recovered  in any abelian symmetric mono\"idal category. This is a bit heavy to write, because we have to express all the above calculus in terms of maps and tensor product of maps, but there is no other obstruction.

So, we have worked above in the symmetric mono\"idal category of $\Bbbk$-modules. For the graded case, it suffices to replace it by the category of graded $\Bbbk$-modules with the signed flip~: 
$x\otimes y \mapsto (-1)^{|x||y|}y\otimes x.$
\end{remarque}

\begin{remarque}
All the construction works if we replace the symmetric algebra $S(L)$ by the tensor algebra $T(L)$, equipped with the shuffle coproduct.
\end{remarque}

\section{The case of rooted trees and Connes-Kreimer Hopf algebra.}

In this section, we  study the case of the pre-Lie algebra of colored rooted trees described in (1.3.). First, we will describe in terms of trees the $\circ$ product of $S(\mathcal{PL}(X))$ defined in section 2. Next, we will give another proof of the result F. Chapoton and M. Livernet \cite{Ch-L}, which states the freeness of $\mathcal{PL}(X)$. Then, we will show that the opposite of our product $*$ is exactly the dual of the coproduct of Connes and Kreimer. For all this section, we fix a set $X$ of color.

\subsection{The $\boldsymbol\circ$ product for the colored rooted trees. }

Let $X_1$, \dots , $X_k$ be trees and consider the monomial $A=X_1\cdots X_k$. Then $A$ is no longer a tree, but we can look at it as the disjoint union of the trees $T_i$, and we often say that it is a forest. 

\vspace{1mm}
Now, let $A$ be a forest, and let $T_1$, \dots, $T_n$ be trees. Then $A\circ T_1\cdots T_n$ is the sum of all the forest resulting from graftings of the trees $T_i$ on the vertices of $A$ : 
\vspace{-1mm}
\begin{displaymath}
A\circ T_1\cdots T_n \, =\,\sum_{(v_1,\cdots,v_n)\in |A|^n} (\cdots ((A\circ_{v_1}T_1)\circ_{v_2}T_2)\cdots )\circ_{v_n}T_n
\end{displaymath}
\vspace{-3mm}
 
For instance, we will need the fact that :
\vspace{-3mm}
\begin{center}
\hspace{1.5cm}\rput(-1.3,0.5){$\bullet \circ T_1\cdots T_n$\;=}\pstree[treemode=U,treesep=5mm,levelsep=1cm]{\Tc*{2pt}}{\Tcircle{$\scriptstyle T_1$}\Tn $\cdots$\Tcircle{$\scriptstyle T_n$}}
\end{center}

\vspace{1mm}
Obviously, all this makes sense when the trees are colored by a set $X$.

\begin{proposition}{\rm \cite{Ch-L}}\quad

\vspace{1mm}
The pre-Lie algebra $\mathcal{PL}(X)$ defined in (1.3.) is the free pre-Lie algebra generated by $X$.
\end{proposition}

\begin{proof} Let $L$ be a pre-Lie algebra, and $f$ a map from $X$ to $L$. Suppose that $\varphi$ is morphism of pre-Lie algebras from $\mathcal{PL}(X)$ to $L$, which extends $f$. As all our constructions of section 2 are functorial, $S(\varphi)$ is a morphism of Hopf algebras from $S(\mathcal{PL}(X))$ to $S(L)$, which commutes with the $\circ$ product defined in (2.6.). Then, for any $X$-colored trees $T_1$, \dots, $T_n$ and any $x\in X$, we have : 
\vspace{-3mm}
\begin{center}
\hspace{-4cm}\rput(-0.4,0.5){$\varphi \Big($}
\rput(3.1,0.5){$\Big)$}
\rput(5.5,0.5){$=\; f(x)\circ \varphi(T_1)\cdots \varphi(T_n)$}
\pstree[treemode=U,treesep=5mm,levelsep=1cm]{\Tcircle{$\scriptscriptstyle x$}}{\Tcircle{$\scriptstyle T_1$} \Tn $\cdots$\Tcircle{$\scriptstyle T_n$}}
\end{center}

This leads us to a definition of $\varphi$ by necessary condition and by induction on the number of vertices. It remains to check that it is a morphism of pre-Lie algebras, which is again done by induction on the number of vertices of the first argument, and follows from the identity (2.7-$v$) :

\rput(2.8,0.5){$\varphi \Big($}
\rput(6.5,0.5){$\circ T_{n+1}\Big)$}
\rput(9.7,0.5){$=\;\varphi\big((x \circ T_1\cdots T_n)\circ T_{n+1}\big)$}
\rput(7.5,-0.8){$\displaystyle =\;\varphi\big((x\circ T_{n+1})\circ T_1\cdots T_n\big ) +  \sum_{i=1}^n \varphi\big(x\circ T_1\cdots (T_i\circ T_{n+1})\cdots T_n\big )\hbox{\rm , by (2.7-}v$),}
\rput(7.5,-1.9){$\displaystyle =\;\varphi(x\circ T_{n+1})\circ \varphi(T_1)\cdots \varphi(T_n) +  \sum_{i=1}^n \varphi(x)\circ \varphi(T_1)\cdots \varphi(T_i\circ T_{n+1})\cdots \varphi(T_n)\hbox{\rm , by definition of }\varphi$,}
\rput(7.5,-3.1){$\displaystyle =\;\big (\varphi(x)\circ \varphi (T_{n+1})\big )\circ \varphi(T_1)\cdots \varphi(T_n) +  \sum_{i=1}^n \varphi(x)\circ \varphi(T_1)\cdots \big (\varphi(T_i)\circ \varphi(T_{n+1})\big )\cdots \varphi(T_n)\hbox{\rm , by induction}$,}
\hspace{3cm} \pstree[treemode=U,treesep=5mm,levelsep=1cm]{\Tcircle{$\scriptscriptstyle x$}}{\Tcircle{$\scriptstyle T_1$} \Tn $\cdots$\Tcircle{$\scriptstyle T_n$}}

\vspace{3.8cm}

\rput(7.5,0){$=\;\big (\varphi(x)\circ \varphi(T_1)\cdots \varphi(T_n)\big )\circ \varphi (T_{n+1})\hbox{\rm , by (2.7-}v$),}
\rput(3.8,-1.4){$=\;\varphi \Big($}
\rput(9.4,-1.4){$\circ T_{n+1}\Big)$, by definition of $\varphi$.}
\rput(5.8,-1.3){\pstree[treemode=U,treesep=5mm,levelsep=1cm]{\Tcircle{$\scriptscriptstyle x$}}{\Tcircle{$\scriptstyle T_1$} \Tn $\cdots$\Tcircle{$\scriptstyle T_n$}}}

\vspace{1.2cm}
\end{proof}
\vspace{0.5cm}

\begin{proposition}\quad

\vspace{1mm}
Let $\mathcal{PL}$ be the pre-Lie algebra of rooted trees. The $*$ product we have defined on $S(\mathcal{PL})$ is the exact dual of the opposite of the Connes-Kreimer coproduct. 
\end{proposition}

\begin{proof} The coproduct of Connes and Keimer is defined inductively by using the following formula :
\vspace{-2mm}
\begin{displaymath}
\Delta \mathcal B_+ \;=\; \mathcal B_+\otimes \eta\,+\, (1\otimes \mathcal B_+)\Delta
\end{displaymath}

Here, $\eta$ stands for the unit of $S(\mathcal{PL})$ and $\mathcal B_+$ is the map which consists in adding a root to a forest : 

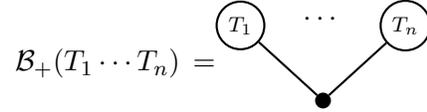
\rput(6.4,0.5){$\mathcal B_+(T_1\cdots T_n)\,= $}
\hspace{7.6cm} \pstree[treemode=U,treesep=5mm,levelsep=1cm]{\Tc*{3pt}}{\Tcircle{$\scriptstyle T_1$} \Tn $\cdots$\Tcircle{$\scriptstyle T_n$}}

Therefore, in order to prove our claim, it suffices to show that our product satisfies the following identity :
\vspace{-2mm}
\begin{displaymath}
\mathcal{B}_-(A*B) \;=\; \varepsilon(A)\mathcal B_-(B)\,+\, \mathcal B_-(A)*B
\end{displaymath}

where $\mathcal B_-$ is the transpose of $\mathcal B^+$. For this, let us first notice that $B_-(A)$ is null for any homogenous element $A$ of $S(\mathcal{PL})$ whose length is different from $1$. Moreover : 
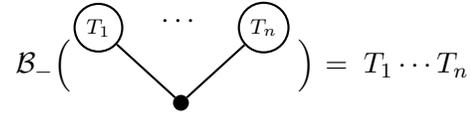
\vspace{-2mm}
\begin{center}
\hspace{-3cm}\rput(0.3,0.5){$\mathcal B_- \Big($}
\rput(3.8,0.5){$\Big)$}
\rput(5,0.5){$=\; T_1\cdots T_n$}
\hspace{7mm}\pstree[treemode=U,treesep=5mm,levelsep=1cm]{\Tc*{3pt}}{\Tcircle{$\scriptstyle T_1$} \Tn $\cdots$\Tcircle{$\scriptstyle T_n$}}
\end{center}

Now, as for any homogenous elements $A$ and $B$ of $S(\mathcal{PL})$, $A\circ B$ has the same length as $A$, $\mathcal{B}_-(A*B)$ is null whenever the length of $A$ is greater than $2$. The above identity is therefore checked in this case, and it remains to get it when $A$ is a constant or a tree. For the constant case, we get the obvious equality $\mathcal B_-(B)=\mathcal B_-(B)$. The following leads us to the conclusion of our proof :
\vspace{-2mm}
\begin{displaymath}
\begin{array}{ccl}
\mathcal B_-\big ((\bullet\circ T_1\cdots T_n)*B\big )&=&\mathcal B_-\big ((\bullet\circ T_1\cdots T_n)\circ B_{(1)}B_{(2)}\big ),\\
&=&\mathcal B_-\big ((\bullet\circ T_1\cdots T_n)\circ B\big ),\\
&=&\mathcal B_-\big (\bullet\circ (T_1\cdots T_n * B)\big )\quad\hbox{\rm by (2.7.{\it v})},\\
&=& T_1\cdots T_n * B
\end{array}
\vspace{-7mm}
\end{displaymath}
\end{proof}

\section{The link with symmetric brace algebras.}

Symmetric brace algebras were introduced by T. Lada and M. Markl in \cite{L-M}. In this part, we will show, as an application of section 2, that symmetric brace algebras and pre-Lie algebras are the same. First, let us recall the definition of symmetric brace algebras and some basic facts.

\begin{definition}\cite{L-M}

A symmetric brace algebra is a vector space $V$ equipped with a brace :
\vspace{-2mm}
\begin{displaymath}
\begin{array}{ccc}
V\otimes S(V)&\longrightarrow & V\\
X\otimes A&\longmapsto & X\{A\}
\end{array}
\end{displaymath}

satisfying the following identity : 
\vspace{-2mm}
\begin{displaymath}
\begin{array}{c}
X\{1\}=X,\\
X\{Y_1,\cdots ,Y_n\}\{A\} = X\{Y_1\{A_{(1)}\},\cdots,Y_n\{A_{(n)}\},A_{(n+1)}\}.
\end{array}
\end{displaymath}
\vspace{-2mm}
\end{definition}

\begin{proposition}{\rm \cite{L-M}}

\vspace{1mm}
Let $V$ be a symmetric brace algebra and define : 
\vspace{-2mm}
\begin{displaymath}
X\circ Y\, =\, X\{Y\}.
\end{displaymath}

\vspace{-2mm}
Then $(V,\circ)$ is a pre-Lie algebra.
\end{proposition}

\begin{proposition}\quad

\vspace{1mm}
Let $V$ be a symmetric brace algebra. Then, the following identity holds :

\vspace{-2mm}
\begin{displaymath}
X\{A\}\, =\, X\circ A.
\end{displaymath}

\vspace{-2mm}
where $\circ$ was defined in (2.4.) and (2.6.).
\end{proposition}

\begin{proof} This is done by induction on the length of $A$. From the identity of (4.1.), we get 
\vspace{-2mm}
\begin{displaymath}
\begin{array}{ccc}
(X\circ A )\circ Y & = & X\circ (A\circ X) + X\{AY\}
\end{array}
\end{displaymath}

\vspace{-2mm}
and we conclude with (2.6.{\it ii})
\end{proof}

\begin{proposition}\quad

\vspace{1mm}
Let $(L,\circ)$ be a pre-Lie algebra. The following braces on $L$ : 
\vspace{-2mm}
\begin{displaymath}
X\{Y_1\cdots Y_n\} \,=\, X\circ Y_1\cdots Y_n 
\end{displaymath}

\vspace{-2mm}
define a symmetric brace algebra structure on $L$.

\end{proposition}

\begin{proof}
This is a straightforward consequence of (2.7.{\it v}) and (2.6.{\it iii}).
\end{proof}

\begin{corollaire}\quad

The categories of symmetric brace algebras and pre-Lie algebras are isomorphic.\qed
\end{corollaire}

\section{Link with brace algebras and $\boldsymbol{\mathcal B}_{\boldsymbol{\infty}}$-algebras.}

In this section, we will extend the construction of part 2 in the brace algebra setting. This will lead us to an already known fact that was used by Tamarkin in his proof of the Kontsevich formality theorem \cite{T} \cite{Hi} : any brace algebra structure on a vector space extends to a $\mathcal B_{\infty}$-algebra one. Then, we will recover the Hopf algebra of planar rooted trees introduced by Foissy \cite{F}

Let us first recall the definition of brace algebra \cite{G-J}. For this, we will use the usual structure of non cocommutative coalgebra on $T(V)$.

\begin{definition}\quad

\vspace{1mm}
Let $V$ be a vector space. A brace algebra structure on $V$ is given by a map :
\vspace{-2mm}
\begin{displaymath}
\begin{array}{ccc}
V\otimes T(V) & \longrightarrow & V\\
X\otimes A & \longmapsto & X\{A\}
\end{array}
\end{displaymath}

\vspace{-2mm}
satisfying the following relations :
\vspace{-2mm}
\begin{displaymath}
\begin{array}{c}
X\{1\}=X,\\
X\{Y_{1} \cdots  Y_{n}\}\{A\} = X\{A_{(1)}\,Y_{1}\{A_{(2)}\} \,A_{(3)}\,\cdots \,A_{(2n-1)}\,Y_{n}\{A_{(2n)}\}\,A_{(2n+1)}\},
\end{array}
\end{displaymath}
where the Sweedler notations are here used for the deconcatenation coproduct of $T(V)$ : 
\vspace{-2mm}
\begin{displaymath}
\Delta (X_{1}\cdots X_{n})=\sum_{i=0}^{n}X_{1}\cdots X_{i}\otimes X_{i+1}\cdots X_{n}.
\end{displaymath}
\end{definition}

Now, let us define an associative product, similar to the one of (2.8.), on the tensor cogebra $T(V)$ of any brace algebra $V$ :

\begin{proposition}\quad

\vspace{1mm}
Let $V$ be a brace algebra. We define by induction the following $*$ product on $T(V)$ : 
\vspace{-2mm}
\begin{displaymath}
\begin{array}{l}
1*B  \,= \, B\\
XA*B \, = \, B_{(1)}\,X\{B_{(2)}\}\,(A*B_{(3)})
\end{array}
\end{displaymath}

\vspace{-1mm}
for $A$ and $B$ in $T(V)$ and $X$ in $V$. This product is associative and makes $T(V)$ into a Hopf algebra.
\end{proposition}

\begin{proof}
First, let us recall a formula relating the deconcatenation coproduct and the concatenation product of $T(V)$ :
\vspace{-2mm}
\begin{displaymath}
\Delta (AXB)\,=\,A_{(1)}\otimes A_{(2)}XB + AXB_{(1)}\otimes B_{(2)}
\end{displaymath}

where $A$ and $B$ are in $T(V)$, and $X$ is in $V$. Now, we can prove the compatibility of the $*$ product with $\Delta$ by induction :
\vspace{-2mm}
\begin{displaymath}
\begin{array}{ccl}
\Delta(XA*B)&=& B_{(1)}\otimes B_{(2)}\,X\{B_{(3)}\}\,(A*B_{(4)}) + B_{(1)}\,X\{B_{(2)}\}\,(A_{(1)}*B_{(3)})\otimes A_{(2)}*B_{(4)}\\
&=&B_{(1)}\otimes (XA*B_{(2)}) + (XA_{(1)}*B_{(1)})\otimes (A_{(2)}*B_{(2)})\\
&=&(1\otimes XA + XA_{(1)}\otimes A_{(2)})* B_{(1)}\otimes B_{(2)}.
\end{array}
\end{displaymath}

Now, let us state a formula relating the $*$ and the concatenation products :
\vspace{-2mm}
\begin{displaymath}
AYB*C\,=\, (A*C_{(1)})\,Y\{C_{(2)}\}\,(B*C_{(3)}).
\end{displaymath}

\vspace{-2mm}
This is done by induction on the length of the monomial $A$ : 
\vspace{-2mm}
\begin{displaymath}
\begin{array}{ccl}
XAYB*C&=& (A*C_{(1)})\,X\{C_{(2)}\}\,(AYB*C_{(3)})\\
&=& C_{(1)}\,X\{C_{(2)}\}\, (A*C_{(3)})\,Y\{C_{(4)}\}\,(A*C_{(5)})\\
&=& (XA * C_{(1)})\,Y\{C_{(2)}\}\,(A*C_{(3)}).
\end{array}
\end{displaymath}

\vspace{-2mm}
Next, we give an interpretation of the brace identity of (5.1.) in terms of $*$ product :
\vspace{-2mm}
\begin{displaymath}
X\{A\}\{B\}\,=\, X\{A*B\}.
\end{displaymath}

\vspace{-2mm}
Finally, we prove the associativity of $*$ by induction on the length of $A$ : 
\vspace{-2mm}
\begin{displaymath}
\begin{array}{ccl}
(XA*B)*C & = & (B_{(1)}\,X\{B_{(2)}\}\,(A*B_{(3)}))*C\\
&=& (B_{(1)}*C_{(1)})\,X\{B_{(2)}\}\{C_{(2)}\}\,((A*B_{(3)})*C_{(3)})\\
&=& (B_{(1)}*C_{(1)})\,X\{B_{(2)}*C_{(2)}\}\,(A*(B_{(3)}*C_{(3)}))\\
&=& XA * (B*C).
\end{array}
\end{displaymath}
\vspace{-2mm}
\end{proof}

\begin{remarque}
The above $*$ product satisfies the following formula :
\vspace{-2mm}
$$
X_1\cdots X_n*B=B_{(1)}\,X_1\{B_{(2)}\}\,B_{(3)}\,\cdots  \,B_{(2n-1)}\,X_n\{B_{(2n)}\}\,B_{(2n+1)}
$$

Therefore, the brace identity can take the following form :
\vspace{-2mm}
$$
X\{B\}\{C\}=X\{B*C\}
$$
\end{remarque}

\begin{remarque}
Up to a shift, we have shown that any brace algebra is naturally a $\mathcal B_{\infty}$-algebra. Let us recall that a $\mathcal B_{\infty}$-algebra structure on $V$ is given by a structure of Hopf  algebra on $T(V)$, where $T(V)$ is equipped with the deconcatenation coproduct. This fact was well-known and used in Tamarkin's proof of Kontsevich formality theorem \cite{T}.  It can be stated by using the cofree universal property of $T(V)$ \cite{G-J}\cite{Hi}. By the same way, we could have used the cofree universal property of $S(V)$ to state (2.9.) starting from (2.5.).
\end{remarque}

\subsection{The link with Connes-Kreimer like Hopf algebras for planar trees. }

In \cite{F}, Foissy introduced a Connes-Kreimer like Hopf algebra with planar rooted trees. He considered the linear space freely generated by the rooted planar trees and constructed on its tensor algebra a copruduct, making this tensor algebra neither commutative, nor cocommutative Hopf algebra. Following \cite{Ch2} and \cite{R}, it turns out that the space of planar rooted is a (free) brace algebra. We propose to see that the Foissy's coproduct is the exact dual of the $*$ product of (5.2.).

\vspace{3mm}
For a vertex $x$ of a rooted tree $T$, we will denote $In(x)$ the set of its upgoing edges. 

\begin{center}
$\pscirclebox[linestyle=dashed,dash=2pt 2pt]{
\overbrace{\vbox{\hsize=1,8cm\pstree[treemode=U,treesep=5mm,levelsep=1cm]{\Tp*}{\pstree{\Tcircle[edge=\dashededge]{$ \scriptscriptstyle x$}}{\Tp*\Tn \dots \Tp*}}}}^{In(x)}}
$
\end{center}

\begin{df}\quad

A rooted planar tree is rooted tree $T$ where, for every vertex $x$ of $T$, $In(x)$ is equipped with a total order $\leq_x$. 

\vspace{1mm}
Obviously, for a given set $X$, an $X$-coloured planar rooted tree is a planar rooted tree $T$ equipped with a \emph{color} map from $|T|$ to $X$.

\end{df}

\begin{prop} {\rm \cite{Ch2}, \cite{R}}

Let $\mathcal{BR}(X)$ be the linear space generated by the $X$-colored planar rooted trees. There is a unique brace algebra structure on $\mathcal{BR}(X)$  such that : 
\vspace{2mm}

\rput(6.2,0.5){$\bullet\{T_1\cdots T_n\}\,= $}
\hspace{7.6cm} \pstree[treemode=U,treesep=5mm,levelsep=1cm]{\Tc*{3pt}}{\Tcircle{$\scriptstyle T_1$} \Tn $\cdots$\Tcircle{$\scriptstyle T_n$}}

\vspace{2mm}
for any element $\bullet$ of $X$ and any $X$-colored planar rooted trees $T_1$, \dots, $T_n$.  \hfill $\square$
\end{prop}

\vspace{3mm}
This brace structure on planar rooted trees was exhibited by Chapoton in \cite{Ch2}. An inductive proof of this proposition can be performed, but we prefer to give a combinatorial description of this braces at the end of the paper.

\begin{thm}\quad

The coproduct of the Hopf algebra of $X$-colored planar rooted trees introduced by Foissy in \cite{F} is dual to the opposite of the $*$ product defined in (5.2.). 

\end{thm}

\begin{proof} It is almost the same as the proof of (3.2.). Indeed, Foissy's coproduct is also defined inductively by using the formulas :
\vspace{-2mm}
\begin{displaymath}
\Delta \mathcal B_x^+ \;=\; \mathcal B_x^+\otimes \eta\,+\, (1\otimes \mathcal B_x^+)\Delta,
\end{displaymath}

where $x$ lies in $X$. But now, $\eta$ stands for the unit of $T(\mathcal{BR})$ and $\mathcal B_x^+$ is the map which consists in adding an $x$-colored root to a forest : 

\rput(6.4,0.5){$\mathcal B_x^+(T_1\cdots T_n)\,= $}
\hspace{7.6cm} \pstree[treemode=U,treesep=5mm,levelsep=1cm]{\Tcircle{$\scriptscriptstyle x$}}{\Tcircle{$\scriptstyle T_1$} \Tn $\cdots$\Tcircle{$\scriptstyle T_n$}}

\vspace{2mm}
Therefore, in order to prove our claim, it suffices to show that the $*$ product satisfies the following identity :
\vspace{-2mm}
\begin{displaymath}
\mathcal{B}_x^-(A*B) \;=\; \varepsilon(A)\mathcal B_x^-(B)\,+\, \mathcal B_x^-(A)*B
\end{displaymath}

where $\mathcal B_x^-$ is the transpose of $\mathcal B_x^+$. For this, let us first notice that $B_x^-(A)$ is null for any homogenous element $A$ of $T(\mathcal{BR})$ whose length is different from $1$. Moreover, for any $x$ and $\bullet$ in $X$ : 
\begin{center}
\hspace{-3cm}\rput(0.3,0.5){$\mathcal B_x^- \Big($}
\rput(3.8,0.5){$\Big)$}
\rput(5.3,0.5){$=\; \delta_{x,\bullet}\, \, T_1\cdots T_n$}
\hspace{7mm}\pstree[treemode=U,treesep=5mm,levelsep=1cm]{\Tc*{3pt}}{\Tcircle{$\scriptstyle T_1$} \Tn $\cdots$\Tcircle{$\scriptstyle T_n$}}
\end{center}

It follows from (5.2.) that, for any homogenous elements $A$ and $B$ of $T(\mathcal{BR})$, the length of any monomial composing  $A*B$ is greater than the length of $A$. Therefore, $\mathcal B_x^-(A*B)$ is null whenever $A$ is of length greater than $2$. The above identity is therefore checked in this case, and it remains to get it when $A$ is a constant or a tree. For the constant case, we get the obvious equality $\mathcal B_x^-(B)=\mathcal B_x^-(B)$. The following leads us to the conclusion of our proof :
\vspace{-2mm}
\begin{displaymath}
\begin{array}{ccl}
\mathcal B_x^-\big (\bullet\{ T_1\cdots T_n\}*B\big )&=&\mathcal B_x^-\big (B_{(1)}\,(\bullet\{ T_1\cdots T_n\}\{B_{(2)}\})\,B_{(3)}\big ),\\
&=&\mathcal B_x^-\big (\bullet\{ T_1\cdots T_n\}\{B\}\big ),\\
&=&\mathcal B_x^-\big (\bullet\{T_1\cdots T_n * B\}\big )\quad\hbox{\rm by (5.3.),}\\
&=& \delta_{x,\bullet}\,\,\, T_1\cdots T_n * B.
\end{array}
\vspace{-7mm}
\end{displaymath}
\end{proof}

\begin{rk}\quad

Therefore, we do recover the Hopf algebras of Foissy and the Hopf algebras of Loday-Ronco and of Frabetti-Brouder, since these last two are isomorphic to Foissy's ones, as proved by Foissy in \cite{F}.
\end{rk}

The end of the paper is now devoted to a combinatorial description of the braces claimed in (5.1.2.)

\subsection{A combinatorial description of the braces of planar rooted trees.}

First, let us define a total order on the edges of a planar rooted tree :

\begin{df}\quad 

Let $T$ be a planar rooted tree. A path of $T$ is a sequence of adjacent edges of $T$, and we will say it is rooted if it starts from the root.

\vspace{1mm}
Let $l= e_1 \cdots e_p$ and $l'=e'_1\cdots e'_q$ be two rooted path. We will say that $l\leq l'$ if $l$ is a subpath of $l'$ or if $l$ goes to the left of $l'$ :
$$\left\{
\begin{array}{l}
p\leq q\;\hbox{ and } e_k=e'_k \;\hbox{ for all } k\leq p,\\
\hbox{or}\\
\hbox{there is an index } i \;\hbox{and a vertex } x \;\hbox{such that } e_i<_x e'_i \;\hbox{in } In(x).
\end{array}
\right.$$

Now, let $e$ and $e'$ be two edges of $T$. Each of them belongs to a unique path starting from the root, which we respectively denote $l$ and $l'$. We will say that $e\leq e'$ whenever $l\leq l'$.
\end{df}

It is easy to check that this defines a total order on the edges of $T$. 

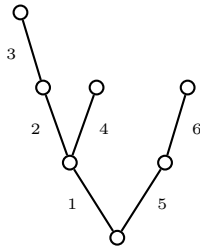
\vspace{-3mm}
\begin{center}
\pstree[treemode=U,treesep=5mm,levelsep=1cm]{\Tc{3pt}}{
 \pstree{\Tc{3pt}\tlput{$\scriptscriptstyle 1$}}{\pstree{\Tc{3pt}\tlput{$\scriptscriptstyle 2$}}{\Tc{3pt}\tlput{$\scriptscriptstyle 3$} \Tn} \Tc{3pt}\trput{$\scriptscriptstyle 4$}}
\pstree{\Tc{3pt}\trput{$\scriptscriptstyle 5$}}{\Tn\Tc{3pt}\trput{$\scriptscriptstyle 6$}}}
\end{center}
\vspace{-2mm}
Now, let us define the sectors of a planar rooted tree :

\vspace{3mm}
\begin{df} \quad

Let $T$ be a planar rooted tree, and $x$ be a vertex of $T$. We define the set $\mathcal S(x)$ of sectors around $x$ to be $In(x)\amalg \{\infty_x\}$ : 
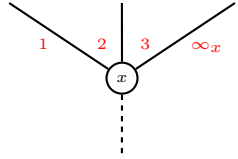
\vspace{-1mm}
\begin{center}
\pstree[treemode=U,treesep=5mm,levelsep=1cm]{\Tp*}{\pstree{\Tcircle[edge=\dashededge]{$ \scriptscriptstyle x$}}{ \Tp*\tlput{$\scriptscriptstyle\red 1$}\Tn\Tn\Tp*\tlput{$\scriptscriptstyle\red 2$}\Tn\Tn\Tp*\tlput{$\scriptscriptstyle\red 3$}\trput{$\scriptscriptstyle\red \infty_x$}}}
\end{center}
\vspace{-2mm}
The total order of $In(x)$ extends naturally to a total order of $In(x)\amalg\mathcal S(x)$. We denote $\mathcal S(T)$ the disjoint union of all $\mathcal S(x)$ for all the vertices $x$ of $T$.
\end{df}

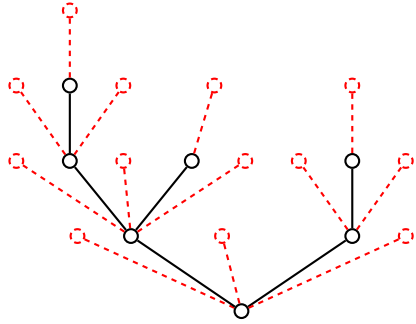
\vspace{-5mm}
\begin{center}
\pstree[treemode=U,treesep=5mm,levelsep=1cm]{\Tc{3pt}}
       {
       \Tc[linestyle=dashed,linecolor=red,edge=\reddashededge]{3pt}
       \pstree{\Tc{3pt}}
              {
              \Tc[linestyle=dashed,linecolor=red,edge=\reddashededge]{3pt}
              \pstree{\Tc{3pt}}
                     {
                     \Tc[linestyle=dashed,linecolor=red,edge=\reddashededge]{3pt}
                     \pstree{\Tc{3pt}}
                            {
                            \Tc[linestyle=dashed,linecolor=red,edge=\reddashededge]{3pt}
                            }
                     \Tc[linestyle=dashed,linecolor=red,edge=\reddashededge]{3pt}
                      } 
              \Tc[linestyle=dashed,linecolor=red,edge=\reddashededge]{3pt}
              \pstree{\Tc{3pt}}
                     {
                     \Tn
                     \Tc[linestyle=dashed,linecolor=red,edge=\reddashededge]{3pt}
                     }
              \Tc[linestyle=dashed,linecolor=red,edge=\reddashededge]{3pt}
              }
       \Tn
       \Tc[linestyle=dashed,linecolor=red,edge=\reddashededge]{3pt}
       \pstree{\Tc{3pt}}
              {
              \Tc[linestyle=dashed,linecolor=red,edge=\reddashededge]{3pt}
              \pstree{\Tc{3pt}}
                     {
                     \Tc[linestyle=dashed,linecolor=red,edge=\reddashededge]{3pt}
                     }
              \Tc[linestyle=dashed,linecolor=red,edge=\reddashededge]{3pt}
              }
       \Tc[linestyle=dashed,linecolor=red,edge=\reddashededge]{3pt}
       }
\end{center}

On the above picture, the red dashed lines represent the elements of $\mathcal S(T)$, where $T$ is the plain black underlying tree. As you can see, we can  associate another rooted tree $\overline{T}$ corresponding to the above dashed red and plain black tree to any planar rooted tree $T$. More precisely, the set of vertices of $\overline{T}$ is $|T|\amalg \mathcal S(T)$ and, for every vertex $x$ of $T$, the set of upgoing edges of $x$ is $In(x)\amalg \mathcal S(x)$ equipped with the order of (5.2.2.).

\vspace{1mm}
\quad Now, by (5.2.1.), we have a total order on the edges of $\overline{T}$. We restrict it in order to get a total order on $\mathcal S(T)$. 
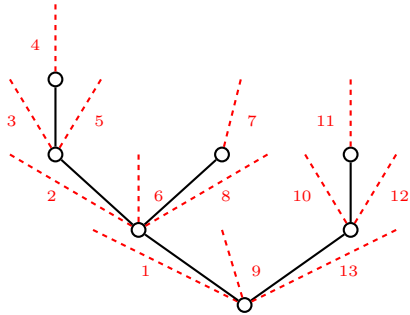
\vspace{-3mm}
\begin{center}
\pstree[treemode=U,treesep=5mm,levelsep=1cm]{\Tc{3pt}}
       {
       \Tp*[edge=\reddashededge]\tlput{$\red\scriptscriptstyle 1$}
       \pstree{\Tc{3pt}}
              {
              \Tp*[edge=\reddashededge]\tlput{$\red\scriptscriptstyle 2$}
              \pstree{\Tc{3pt}}
                     {
                     \Tp*[edge=\reddashededge]\tlput{$\red\scriptscriptstyle 3$}
                     \pstree{\Tc{3pt}}
                            {
                            \Tp*[edge=\reddashededge]\tlput{$\red\scriptscriptstyle 4$}
                            }
                     \Tp*[edge=\reddashededge]\trput{$\red\scriptscriptstyle 5$}
                      } 
              \Tn
              \Tp*[edge=\reddashededge] \trput{$\red\scriptscriptstyle 6$}
              \Tn
              \pstree{\Tc{3pt}}
                     {
                     \Tn
                     \Tp*[edge=\reddashededge]\trput{$\red\scriptscriptstyle 7$}
                     }
              \Tp*[edge=\reddashededge]\trput{$\red\scriptscriptstyle 8$}
              }
       \Tn
       \Tp*[edge=\reddashededge]\trput{$\red\scriptscriptstyle 9$}
       \pstree{\Tc{3pt}}
              {
              \Tp*[edge=\reddashededge]\tlput{$\red\scriptscriptstyle 10$}
              \pstree{\Tc{3pt}}
                     {
                     \Tp*[edge=\reddashededge]\tlput{$\red\scriptscriptstyle 11$}
                     }
              \Tp*[edge=\reddashededge]\trput{$\red\scriptscriptstyle 12$}
              }
       \Tp*[edge=\reddashededge]\trput{$\red\scriptscriptstyle 13$}
       }
\end{center}
  
\begin{df} Grafting of a planar forest along a sector.

Let $T$, $T_1$, \dots, $T_k$ be $X$-colored planar rooted trees and $s\in \mathcal S(x)$ be a sector of $T$.  The grafting of $T_1\cdots T_n$ on $T$ along $s$ is the planar rooted tree denoted $T\circ_s T_1\cdots T_n$ defined as follow :

\vspace{1mm}
\qquad\qquad\quad (i) $|T\circ_s T_1\cdots T_k|=|T|\amalg |T_1|\amalg \cdots \amalg |T_k|$,

\vspace{1mm}
\qquad\qquad\quad (ii)  $T$ and all the $T_i$ are planar subtrees of $T\circ_s T_1\cdots T_k$,

\vspace{1mm}
\qquad \qquad\quad (iii) for all $i$, an  edge $e_i$  goes up from $x$ to the root $T_i$,

\vspace{1mm}
\qquad \qquad\quad (iv) for every $e$ in $In(x)$, $e<e_i$ whenever $e<s$ and $e_i<e$ whenever $s<e$.
\end{df}

When $T$ is the above planar rooted tree and $s$ is its sixth sector, then $T\circ_s T_1\cdots T_k$ is the following tree :
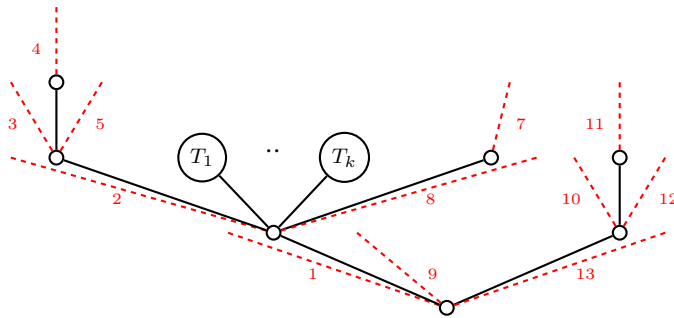
\vspace{-3mm}
\begin{center}
\pstree[treemode=U,treesep=5mm,levelsep=1cm]{\Tc{3pt}}
       {
       \Tp*[edge=\reddashededge]\tlput{$\red\scriptscriptstyle 1$}
       \pstree{\Tc{3pt}}
              {
              \Tp*[edge=\reddashededge]\tlput{$\red\scriptscriptstyle 2$}
              \pstree{\Tc{3pt}}
                     {
                     \Tp*[edge=\reddashededge]\tlput{$\red\scriptscriptstyle 3$}
                     \pstree{\Tc{3pt}}
                            {
                            \Tp*[edge=\reddashededge]\tlput{$\red\scriptscriptstyle 4$}
                            }
                     \Tp*[edge=\reddashededge]\trput{$\red\scriptscriptstyle 5$}
                      } 
              \Tn
              \Tn
              \Tcircle{$\scriptstyle T_1$} 
              \Tn 
              $\cdot\cdot$ 
              \Tcircle{$\scriptstyle T_k$}
              \Tn
              \Tn
              \pstree{\Tc{3pt}}
                     {
                     \Tn
                     \Tp*[edge=\reddashededge]\trput{$\red\scriptscriptstyle 7$}
                     }
              \Tp*[edge=\reddashededge]\trput{$\red\scriptscriptstyle 8$}
              }
       \Tn
       \Tp*[edge=\reddashededge]\trput{$\red\scriptscriptstyle 9$}
       \pstree{\Tc{3pt}}
              {
              \Tp*[edge=\reddashededge]\tlput{$\red\scriptscriptstyle 10$}
              \pstree{\Tc{3pt}}
                     {
                     \Tp*[edge=\reddashededge]\tlput{$\red\scriptscriptstyle 11$}
                     }
              \Tp*[edge=\reddashededge]\trput{$\red\scriptscriptstyle 12$}
              }
       \Tp*[edge=\reddashededge]\trput{$\red\scriptscriptstyle 13$}
       }
\end{center}

Now, let $T_1,\cdots ,T_n$ be colored planar rooted trees and let $f$ be an increasing map from $\{1,\cdots ,n\}$ to $\mathcal S(T)$. We denote $\{s_1<\cdots< s_m\}$ the image of $f$. We can use $f$ to graft the $T_i$ trees on $T$ as follows : 
\vspace{-2mm}
$$T\circ_f T_1\cdots T_n= \bigg(\cdots \Big(\big(T\circ_{s_1} T_{f^{-1}(s_1)}\big)\circ_{s_2}T_{f^{-1}(s_2)}\Big )\cdots  \bigg ) \circ_{s_m}T_{f^{-1}(s_m)},$$

\vspace{-1mm}
where $T_{f^{-1}(s_k)}=T_i\cdots T_{j}$ if $f^{-1}(s_k)=\{i<\cdots <j\}$.

\vspace{3mm}
\begin{df} {\rm The braces of} $\mathcal{BR}(X)$.

Let $T$, $T_1,\cdots,T_n$ be $X$-colored planar rooted trees. We define :
\vspace{-2mm}
$$T\{T_1\cdots T_n\}=\sum_{(f)} T\circ_f T_1\cdots T_n$$

\vspace{-3mm}
where the sum is taken over all increasing maps $f$ from $\{1,\cdots ,n\}$ to $\mathcal S(T)$.
\end{df}


\end{document}